\documentclass[11pt]{article}
\usepackage{amsmath,amssymb,amsfonts,diagrams}
\newenvironment{keywords}{\noindent\small {\it Keywords\/}:}{\vskip 4pt}
\newenvironment{classification}{\noindent\small 2000 {\it Mathematics Subject
Classification\/}:}{\vskip 12pt}

\newcommand{\posints}{{\mathbb N}}

\newcommand{\tensor}{\otimes}

\newcommand{\Tensor}{\hat{\otimes}}

\newcommand{\wTensor}{\check{\otimes}}

\newcommand{\cstar}{{C^\ast}}

\newcommand{\id}{{\mathrm{id}}}

\newcommand{\A}{{\mathfrak A}}
\newcommand{\B}{{\mathfrak B}}
\newcommand{\Hilbert}{{\mathfrak H}}
\newcommand{\Kilbert}{{\mathfrak K}}
\newcommand{\M}{{\mathfrak M}}

\newcommand{\CB}{{\cal CB}}

\newcommand{\VN}{\operatorname{VN}}

\renewcommand{\baselinestretch}{1.2}
\newcommand{\dated}{\mbox{} \hfill {\small [{\tt \today}]}} 
\usepackage{amsthm,enumerate}
\theoremstyle{plain}
\newtheorem{theorem}{Theorem}
\newtheorem*{proposition}{Proposition}
\newtheorem*{corollary}{Corollary}
\newtheorem*{lemma}{Lemma}
\theoremstyle{definition}
\newtheorem{definition}{Definition}
\theoremstyle{remark}
\newtheorem{example}{Example}
\newcommand{\WAP}{\mathrm{WAP}}
\title{Co-representations of Hopf--von Neumann algebras \\ on operator spaces other than column Hilbert space}
\author{\textit{Volker Runde}}
\date{}
\begin{document}
\maketitle
\begin{abstract}
Recently, M.\ Daws introduced a notion of co-representation of abelian Hopf--von Neumann algebras on general reflexive Banach spaces. In this note, we show that this notion cannot be extended beyond subhomogeneous Hopf--von Neumann algebras. The key is our observation that, for a von Neumann algebra $\M$ and a reflexive operator space $E$, the normal spatial tensor product $\M \bar{\tensor} \CB(E)$ is a Banach algebra if and only if $\M$ is subhomogeneous or $E$ is completely isomorphic to column Hilbert space.
\end{abstract}
\section*{Introduction}
If $\A$ is a Banach algebra, then its dual space $\A^\ast$ is a Banach $\A$-bimodule through
\[
  \langle x, a \cdot \phi \rangle := \langle xa, \phi \rangle \quad\text{and}\quad
  \langle x, \phi \cdot a \rangle := \langle ax, \phi \rangle \qquad (a,x \in \A, \, \phi \in \A^\ast).
\]
A functional $\phi \in \A^\ast$ is said to be \emph{weakly almost periodic} if $\{ a \cdot \phi : a \in \A, \, \| a \| \leq 1 \}$ is relatively weakly compact in $\A^\ast$. There appears to be some asymmetry in the definition of a weakly almost periodic functional, but thanks to Grothendieck's double limit criterion (\cite[Proposition 4]{Gro}), $\phi \in \A^\ast$ is weakly almost periodic if and only if $\{ \phi \cdot a : a \in \A, \, \| a \| \leq 1 \}$ is relatively weakly compact in $\A^\ast$. The collection of all weakly almost periodic functionals on $\A$ is a closed subspace of $\A^\ast$, which we denote by $\WAP(\A)$.
\par 
Let $G$ be a locally compact group, and let $\A$ be its group algebra $L^1(G)$. In this case, it is not difficult to see that $\WAP(\A)$ is nothing but $\WAP(G)$, the commutative $\cstar$-algebra of all weakly almost periodic functions on $G$ (see \cite{Bur} for the definition and properties of $\WAP(G)$). Now, let $\A$ be Eymard's Fourier algebra $A(G)$ with dual $\VN(G)$ (see \cite{Eym}); in this case, we denote $\WAP(\A)$ by $\WAP(\hat{G})$. If $G$ is abelian, then $\WAP(\hat{G})$ is indeed just the weakly almost periodic functions on the dual group $\hat{G}$ and therefore, in particular, is a $\cstar$-subalgebra of $L^\infty(\hat{G}) \cong \VN(G)$. With a little more effort, one can show that $\WAP(\hat{G})$ is a $\cstar$-subalgebra of $\VN(G)$ whenever $G$ has an abelian subgroup of finite index. For general $G$, however, it has been an open question for decades whether or not $\WAP(\hat{G})$ is always a $\cstar$-subalgebra of $\VN(G)$.
\section*{Hopf--von Neumann algebras and co-representations}
Recently, M.\ Daws showed in \cite{Daw} that $\WAP(M(G))$ is a $\cstar$-subalgebra of ${\cal C}_0(G)^{\ast\ast}$ for any $G$, where $M(G)$ is the measure algebra of $G$. In fact, Daws proved a much more general result about abelian Hopf--von Neumann algebras.
\begin{definition}
A \emph{Hopf--von Neumann algebra} is a pair $(\M, \Gamma)$, where $\M$ is a von Neumann algebra, and $\Gamma$ is a \emph{co-multiplication}: a unital, injective, normal $^\ast$-homomorphism $\Gamma \!: \M \to \M \bar{\tensor} \M$ 
which is co-associative, i.e., 
\[
  ( \id \tensor \Gamma) \circ \Gamma = (\Gamma \tensor \id) \circ \Gamma.
\]
\end{definition}
\par 
We call a Hopf--von Neumann algebra $(\M,\Gamma)$ abelian, semidiscrete, etc., if the underlying von Neumann algebra $\M$ has the corresponding property.
\begin{example}
Let $G$ be a locally compact group. 
\begin{enumerate}[(a)]
\item A co-multiplication $\Gamma \!: L^\infty(G) \to L^\infty(G) \bar{\tensor} L^\infty(G) \cong L^\infty(G \times G)$ is given by
\[
  (\Gamma \phi)(x,y) := \phi(xy) \qquad (\phi \in L^\infty(G), \, x,y \in G).
\]
(Restricting $\Gamma$ to ${\cal C}_0(G)$ and then taking section adjoints, yields another co-mul\-ti\-pli\-ca\-tion $\tilde{\Gamma} \!: {\cal C}_0(G)^{\ast\ast} \to {\cal C}_0(G)^{\ast\ast} \bar{\tensor} {\cal C}_0(G)^{\ast\ast}$.) 
\item Let $\lambda \!: G \to {\cal B}(L^2(G))$ be the left regular representation of $G$. Then a co-mul\-ti\-pli\-ca\-tion $\hat{\Gamma} \!: \VN(G) \to \VN(G) \bar{\tensor} \VN(G)$ is given by
\[
  \hat{\Gamma}(\lambda(x)) = \lambda(x) \tensor \lambda(x) \qquad (x \in G).
\]
\end{enumerate}
\end{example}
\par 
Whenever $(\M,\Gamma)$ is a Hopf--von Neumann algebra, the unique predual $\M_\ast$ of $\M$ becomes a Banach algebra via
\[
  \langle x, f \ast g \rangle := \langle \Gamma x, f \tensor g \rangle \qquad (f,g \in \M_\ast, \, x \in \M).
\]
\begin{example}
If $G$ is a locally compact group, then $\ast$ defined in this manner for $(L^\infty(G),\Gamma)$ is the usual convolution product on $L^1(G)$ whereas $\ast$ for $(\VN(G),\hat{\Gamma})$ is the pointwise product on $A(G)$. 
\end{example}
\par 
Any von Neumann algebra $\M$ is a (concrete) operator space, so that $\M_\ast$ is an abstract operator space. (For background on operator space theory, we refer to \cite{ER}, the notations of which we adopt.) If $(\M,\Gamma)$ is a Hopf--von Neumann algebra, then $\Gamma$ is a \emph{complete} isometry. Consequently, $(\M_\ast,\ast)$ is not only a Banach algebra, but a \emph{completely contractive Banach algebra} (\cite[p.\ 308]{ER}).
\par 
The main result of \cite{Daw} is:
\begin{theorem} \label{matt}
Let $(\M,\Gamma)$ be an abelian Hopf--von Neumann algebra. Then $\WAP(\M_\ast)$ is a $\cstar$-subalgebra of $\M$. 
\end{theorem}
\par 
At the heart of Daws' proof is the notion of a co-representation of a Hopf--von Neumann algebra. Usually, one considers co-representation on Hilbert spaces:
\begin{definition} \label{corepdef}
Let $(\M,\Gamma)$ be a Hopf--von Neumann algebra. A \emph{co-representation} of $(\M,\Gamma)$ on a Hilbert space $\Hilbert$ is an operator $U \in \M \bar{\tensor} {\cal B}(\Hilbert)$ such that
\begin{equation} \tag{\mbox{$\ast$}} \label{corep}
  (\Gamma \tensor \id)(U) = U_{1,3} U_{2,3}.
\end{equation}
\end{definition}
\par 
Here, $U_{1,3}$ and $U_{2,3}$ are what is called \emph{leg notation}: if $\M$ acts on a Hilbert space $\Kilbert$, then $U_{1,3}$ is the linear operator on the Hilbert space tensor product $\Kilbert \tensor_2 \Kilbert \tensor_2 \Hilbert$ that acts as $U$ on the first and the third factor and as the identity on the second one; $U_{2,3}$ is defined similarly.
\par 
Commonly, co-representations are also required to be unitary, but we won't need that.
\par 
By \cite[Corollary 7.1.5 and Theorem 7.2.4]{ER}, we have the completely isometric identifications
\[
  \M \bar{\tensor} {\cal B}(\Hilbert) = (\M_\ast \Tensor {\cal B}(\Hilbert)_\ast)^\ast = \CB(\M_\ast, {\cal B}(\Hilbert)).
\]
Given an operator $U \in \M \bar{\tensor} {\cal B}(\Hilbert)$, the corresponding map in $\CB(\M_\ast,{\cal B}(\Hilbert))$ is
\begin{equation} \tag{\mbox{$\ast\ast$}} \label{map}
  \M_\ast \to {\cal B}(\Hilbert), \quad f \mapsto (f \tensor \id)(U),
\end{equation}
and (\ref{corep}) is equivalent to (\ref{map}) being a multiplicative map from $(\M_\ast, \ast)$ into ${\cal B}(\Hilbert)$. The advantage of looking at elements of $\M \bar{\tensor} {\cal B}(\Hilbert)$ instead of $\CB(\M_\ast, {\cal B}(\Hilbert))$ is that $\M \bar{\tensor} {\cal B}(\Hilbert)$ is again a von Neumann algebra, so that it makes sense to multiply its elements.
\par 
Let $(\M,\Gamma)$ be an abelian Hopf--von Neumann algebra. Then $\M$ is of the form $L^\infty(X)$ for some measure space $X$. The proof of Theorem \ref{matt} in \cite{Daw} has three main ingredients:
\begin{enumerate}
\item Elements of $\WAP(L^1(X))$ arise as coefficients of representations of $(L^1(X),\ast)$ on reflexive Banach spaces.
\item For a reflexive Banach space $E$, the weak$^\ast$ closure of $L^\infty(X) \tensor {\cal B}(E)$ in ${\cal B}(L^2(X,E))$---denoted by $L^\infty(X) \bar{\tensor} {\cal B}(E)$---can be identified with ${\cal B}(L^1(X), {\cal B}(E))$ (\cite[Proposition 3.2]{Daw}). This identification then allows to define co-re\-pre\-sen\-ta\-tions of $(L^\infty(X),\Gamma)$ on $E$ in analogy with Definition \ref{corepdef}.
\item The product in $L^\infty(X) \bar{\tensor} {\cal B}(E)$ corresponds to the product in $\WAP(L^1(X))$.
\end{enumerate}
Is it possible to adapt this approach to more general Hopf--von Neumann algebras?
\section*{Co-representations on operator spaces}
In \cite{Daw}, Daws uses the symbol $L^\infty(X) \bar{\tensor} {\cal B}(E)$ for the closure of $L^\infty(X) \tensor {\cal B}(E)$ in ${\cal B}(L^2(X,E))$. In operator space theory, the symbol $\bar{\tensor}$ is usually reserved for the \emph{normal spatial tensor product} of dual operator spaces (\cite[p.\ 134]{ER}). Both usages are consistent: for a reflexive Banach space $E$, we have the isometric identifications
\[
  \begin{split}
  {\cal B}(L^1(X), {\cal B}(E)) & = \CB(L^1(X), \CB(\max E)) \\
  & = (L^1(X) \Tensor (\max E \Tensor \min E^\ast))^\ast \\
  & = L^\infty(X) \bar{\tensor} \CB(\max E).
  \end{split}
\]
From \cite[Proposition 3.2]{Daw}, we thus conclude that the product on $L^\infty(X) \tensor {\cal B}(E)$ extends to $L^\infty(X) \bar{\tensor} \CB(\max E)$, turning it into a Banach algebra. More generally, $L^\infty(X) \bar{\tensor} \CB(E)$ is a Banach algebra for any reflexive, homogeneous opeator space $E$, i.e., satisfying $\CB(E) = {\cal B}(E)$ with identical norms.
\par 
Let $\M$ be a semidiscrete von Neumann algebra, and let $E$ be a reflexive operator space. Then we have the completely isometric identifications
\[
  \CB(\M_\ast, \CB(E)) = (\M_\ast \Tensor (E \Tensor E^\ast))^\ast = \M \bar{\tensor} \CB(E).
\]
(We need the semidiscreteness of $\M$ for the second equality: without it, we would not get $\M \bar{\tensor} \CB(E)$, but the \emph{normal Fubini tensor product} $\M \bar{\tensor}_{\cal F} \CB(E)$; see \cite{Kra}.) This suggest that it might be possible to define a notion of co-representation of semidiscrete Hopf--von Neumann algebras on general reflexive operator spaces. Just to meaningfully state the right hand side (\ref{corep}) for some $U \in \M \bar{\tensor} \CB(E)$, we need that $\M \bar{\tensor} \M \bar{\tensor} \CB(E)$ is multiplicatively closed, and to adapt the proof of Theorem \ref{matt} to general (semidiscrete) Hopf--von Neumann algebras, we need that $\M \bar{\tensor} \CB(E)$ is also multiplicatively closed. It all comes down to the question whether or not $\M \bar{\tensor} \CB(E)$ is a Banach algebra for a (semidiscrete) von Neumann algebra and a reflexive operator space $E$.
\par 
For abelian $(\M,\Gamma)$, this is indeed the case, and it is not difficult to extend \cite[Proposition 3.2]{Daw} to a general operator space setting:
\begin{proposition}
Let $L^\infty(X)$ be an abelian von Neumann algebra, and let $E$ be a reflexive operator space. Then the closure of $L^\infty(X) \tensor \CB(E)$ in $\CB(L^2(X,E))$ is isometrically isomorphic to $\CB(L^1(X), \CB(E))$. In particular, the product of $L^\infty(X) \tensor \CB(E)$ extends to $L^\infty(X) \bar{\tensor} \CB(E)$, turning it into a Banach algebra.
\end{proposition}
\par 
Here, the operator space structure on $L^2(X,E)$ is the one obtained through complex interpolation betwenn $L^\infty(X) \wTensor E$ and $L^1(X) \Tensor E$ as described in \cite{Pis}.
\par 
We can even go a little bit beyond the abelian framework. If $\M$ is a subhomogeneous von Neumann algebra, it is of the form
\[
  \M \cong M_{n_1}(\M_1) \oplus \cdots \oplus M_{n_k}(\M_k)
\]
with $n_1, \ldots, n_k \in \posints$ and abelian von Neumann algebras $\M_1, \ldots, \M_k$. This yields:
\begin{corollary}
Let $\M$ be a subhomogeneous von Neumann algebra, and let $E$ be a reflexive operator space. Then the product of $\M \tensor \CB(E)$ extends to $\M \bar{\tensor} \CB(E)$, turning it into a Banach algebra (with bounded, but not necessarily contractive multiplication).
\end{corollary}
\par 
So, if $(\M,\Gamma)$ is a subhomogeneous Hopf--von Neumann algebra, we can meaningfully speak of its co-representations on reflexive operator spaces. But what if we go beyond subhomogeneous von Neumann algebras? For certain operator spaces, this is no problem. Let $\Hilbert$ be a Hilbert space, and let $\Hilbert_c$ be column Hilbert space (\cite[3.4]{ER}). Then $\CB(\Hilbert_c) = {\cal B}(\Hilbert)$ as operator spaces (\cite[Theorem 3.4.1]{ER}), so that 
\[
  \M \bar{\tensor} \CB(\Hilbert_c) = \M \bar{\tensor} {\cal B}(\Hilbert)
\]
is a von Neumann algebra, and co-representations on $\Hilbert_c$ are nothing but co-re\-pre\-sen\-ta\-tions on $\Hilbert$ in the sense of Definition \ref{corepdef}.
\par
As we shall see, this is about as far as we can get:
\begin{theorem} \label{main}
Let $(\M,\Gamma)$ be a Hopf--von Neumann algebra, and let $E$ be a reflexive operator space. Then the following are equivalent:
\begin{enumerate}[\rm (i)]
\item The product of $\M \tensor \CB(E)$ extends to $\M \bar{\tensor} \CB(E)$ turning it into a Banach algebra with bounded, but not necessarily contractive multiplication.
\item The product of $\M \tensor \CB(E)$ extends to $\M \wTensor \CB(E)$ turning it into a Banach algebra with bounded, but not necessarily contractive multiplication.
\item $\M$ is subhomogeneous or $E$ is completely isomorphic to $\Hilbert_c$ for some Hilbert space $\Hilbert$.
\end{enumerate}
\end{theorem}
\par 
For the proof, we require a lemma:
\begin{lemma} 
Let $\A$ and $\B$ be completely contractive Banach algebras such that that $\A$ contains the full matrix algebra $M_n$ as a subalgebra for each $n \in \posints$, and suppose that the product of $\A \tensor \B$ extends to $\A \wTensor \B$, turning it into a Banach algebra with bounded, but not necessarily contractive multiplication. Then $\B$ is completely isomorphic to an operator algebra.
\end{lemma}
\par 
Here, an operator algebra is a closed, but not necessarily self-adjoint subalgebra of ${\cal B}(\Hilbert)$ for some Hilbert space $\Hilbert$.
\begin{proof}
Let $C \geq 1$ be a bound for the multiplication in $\A \wTensor \B$, and note that, for $n \in \posints$, we have a complete isometry $M_n(\B) \cong M_n \wTensor \B$ (\cite[Corollary 8.1.3]{ER}). Since $M_n$ is a subalgebra of $\A$, this means that formal matrix multiplication from $M_n(\B) \times M_n(\B)$ to $M_n(\B)$ is bounded by $C$ for each $n$. By the definition of the \emph{Haagerup tensor product} (see \cite[Chapter 9]{ER}) $\B \tensor^h \B$, this means that multiplication $m \!: \B \tensor \B \to \B$ extends to a completely bounded map $m \!: \B \tensor^h \B \to \B$. Hence, $\B$ is completely isomorphic to an operator algebra by \cite[Theorem 5.2.1]{BLeM}.
\end{proof}
\begin{proof}[Proof of Theorem \emph{\ref{main}}]
We have previously seen that (iii) $\Longrightarrow$ (i) holds, and (i) $\Longrightarrow$ (ii) is obvious.
\par 
(ii) $\Longrightarrow$ (iii): Assume that $\M$ is \emph{not} subhomogeneous. Then the structure theory of von Neumann algebras yields that $\M$ contains $M_n$ as a subalgebra for each $n$. By the Lemma, $\CB(E)$ is thus completely isomorphic to an operator algebra. By \cite[Theorem 3.4]{Ble}, this means that $E$ is completely isomorphic to $\Hilbert_c$ for some Hilbert space $\Hilbert$.
\end{proof}
\renewcommand{\baselinestretch}{1.0}
\renewcommand{\baselinestretch}{1.2}
\dated
\vfill
\begin{tabbing}
\textit{Author's address}: \= Department of Mathematical and Statistical Sciences \\
\> University of Alberta \\
\> Edmonton, Alberta \\
\> Canada T6G 2G1 \\[\medskipamount]
\textit{E-mail}: \> \texttt{vrunde@ualberta.ca} \\[\medskipamount]
\textit{URL}: \> \texttt{http://www.math.ualberta.ca/$^\sim$runde/}
\end{tabbing}           

\end{document}